\documentclass{article}
\usepackage{spconf}
\usepackage{amsmath,graphicx,amssymb,amsthm}
\usepackage{algorithm, algorithmic,caption,subcaption,color}

\newcommand{\C}{\mathbb{C}}

\theoremstyle{remark}

\theoremstyle{example}

\theoremstyle{definition}

\title{A FAST RANDOMIZED KACZMARZ ALGORITHM FOR SPARSE SOLUTIONS OF CONSISTENT LINEAR SYSTEMS}
\name{Hassan Mansour
 and {\"O}zg{\"u}r Y{\i}lmaz \vspace{-0.1in}\thanks{hassanm@cs.ubc.ca; 
oyilmaz@math.ubc.ca}
\thanks{
Both authors were supported in part by the Natural
Sciences and Engineering Research Council of Canada (NSERC)
Collaborative Research and Development Grant DNOISE II
(375142-08). {\"O}. Y{\i}lmaz was also supported in part by an NSERC
Discovery Grant.}}
\address{\vspace{-0.1in}University of British Columbia, Vancouver - BC, Canada}

\begin{document}
%
\maketitle
\begin{abstract}
The Kaczmarz algorithm is a popular solver for overdetermined linear systems due to its simplicity and speed. In this paper, we propose a modification that speeds up the convergence of the randomized Kaczmarz algorithm for systems of linear equations with sparse solutions. The speedup is achieved by projecting every iterate onto a weighted row of the linear system while maintaining the random row selection criteria of Strohmer and Vershynin. The weights are chosen to attenuate the contribution of row elements that lie outside of the estimated support of the sparse solution. While the Kaczmarz algorithm and its variants can only find solutions to overdetermined linear systems, our algorithm surprisingly succeeds in finding sparse solutions to underdetermined linear systems as well. We present empirical studies which demonstrate the acceleration in convergence to the sparse solution using this modified approach in the overdetermined case. We also demonstrate the sparse recovery capabilities of our approach in the underdetermined case and compare the performance with that of $\ell_1$ minimization.
\end{abstract}
\section{Introduction}
\label{sec:intro}

The Kaczmarz method \cite{ka37} is a popular algorithm for finding the solution of an overdetermined system of  linear equations 
\begin{equation}\label{eq:LinSys}
	Ax = b,
\end{equation}
where $A \in \C^{m \times n}$ and $b \in \C^m$. When $A$ is large, it often becomes costly to solve for $x$ by inverting the matrix. The Kaczmarz algorithm iteratively cycles through the rows of $A$, orthogonally projecting the iterate onto the solution space given by each row. Due to its speed, simplicity, and low memory usage, the Kaczmarz method has risen to prominence in applications such as computer tomography and signal and image processing \cite{Gordon_etal:1970,SezanStark:87,natterer2001} 

Recently, a randomized version of the Kaczmarz algorithm was proposed in \cite{StrohmerKacz:2009} that enjoys an expected linear rate of convergence. Denote the rows of A by $a_1,\dots a_m$ and let $b = \left(b(1)\dots b(m)\right)^T$. Contrary to the original Kaczmarz method in \cite{ka37} which sequentially cycles through the rows of $A$, the randomized Kaczmarz (RK) algorithm of Strohmer and Vershynin \cite{StrohmerKacz:2009} randomizes the row selection criteria of the Kaczmarz method by choosing the next row $a_i$ indexed by $i \in \{1,\dots m\}$ with probability $\frac{\|a_i\|_2^2}{\|A\|_F^2}$. A full description of the RK algorithm is shown in Algorithm \ref{alg:RK}. At every iteration $j$, the previous estimate $x_{j-1}$ is orthogonally projected onto the space of all points $u \in \C^n$ defined by the hyperplane $\langle a_i, u\rangle = b(i)$.

\begin{algorithm}\caption{Randomized Kaczmarz algorithm \cite{StrohmerKacz:2009}}
\begin{algorithmic}[1]\label{alg:RK}
\STATE \textbf{Input} $b = Ax$, the $m \times n$ matrix $A$, maximum iterations $J$
\STATE \textbf{Output} $x_j$
\STATE \textbf{Initialize} $j = 0$, $x_0 = 0$  
\WHILE{$j \leq J$}
\STATE $j = j + 1$
\STATE Choose the row vector $a_i$ indexed by $i \in \{1,\dots m\}$ with probability $\frac{\|a_i\|_2^2}{\|A\|_F^2}$
\STATE {\large$x_j = x_{j-1} + \frac{b(i) - \langle a_i, x_{j-1}\rangle}{\|a_i\|_2^2} a_i^T$}
\ENDWHILE
\end{algorithmic}
\end{algorithm}

In many applications in tomography and image processing, the solution vector $x$ is sparse or admits a sparse representation in some transform domain. In this paper, we propose a modification of the RK algorithm that accelerates the convergence of RK when the solution to $Ax = b$ is sparse with support $S_0$. Our sparse randomized Kaczmarz (SRK) algorithm uses the same row selection criteria as RK, but achieves faster convergence by projecting onto a weighted row $a_i$ of the matrix $A$. The weighting is based on identifying a support estimate $S$ for $x$ and gradually scaling down the entries of $a_i$ that lie outside of $S$ by a factor equal to $1/\sqrt{j}$. As the number of iterations becomes large, the weight $\frac{1}{\sqrt{j}} \rightarrow 0$ and our algorithm begins to resemble the RK algorithm applied to the reduced system
\begin{equation}\label{eq:RedLinSys}
	A_Sx_S = b.
\end{equation}
So long as the support estimate $S$ remains a superset of the true support $S_0$, the solution to the reduced system \eqref{eq:RedLinSys} will be identical to the solution of the full system \eqref{eq:LinSys}. Moreover, Strohmer and Vershynin showed that the convergence rate of the RK algorithm is inversely proportional to the scaled condition number of the matrix $A$, $\kappa(A)$. More recently, Chen and Powell \cite{ChenPowell:12} also proved almost sure linear convergence of the RK algorithm. Consequently, under the assumption that $S_0 \subseteq S$, our algorithm converges at a rate inversely proportional to $\kappa(A_S) \leq \kappa(A)$\footnote{The condition number of a submatrix is always smaller than or equal to the condition number of the full matrix \cite{Thompson1975}.}, which explains the acceleration in the convergence rate of SRK. A full proof of the convergence of our algorithm is still underway and we leave detailed explanations to a future paper.

In Section \ref{sec:PriorWork}, we discuss how our work relates to previous work on accelerating the randomized Kaczmarz algorithm. We then present in Section \ref{sec:SparseRandKacz} the details of our sparse randomized Kaczmarz algorithm. In Section \ref{sec:NumExp}, we provide numerical experiments that demonstrate the acceleration in convergence rate to the sparse solution of overdetermined linear systems. We also demonstrate in the same section the sparse recovery capabilities of SRK in the case where the linear system is underdetermined. 

\textbf{Remark:} A slight modification to our algorithm allows it to effectively handle inconsistent systems as well as compressible signals $x$. However, we omit a discussion of these cases due to space limitations.\vspace{-0.2in}

\section{Relation to prior work}\label{sec:PriorWork}
Several works in the literature have proposed ways of accelerating the convergence of the randomized Kaczmarz algorithm for overdetermined linear systems. In \cite{elne11}, acceleration is achieved by projecting the rows of $A$ onto a lower dimensional space and then performing the Kaczmarz iteration with respect to these low dimensional vectors. Another approach proposed in \cite{NeedellWard:12,NeedellTropp:12} involves randomly selection blocks of rows of $A$ in the Kaczmarz projection. Note that our weighting approach can be combined with any of the above mentioned works to speed up their convergence rate when the desired solution is sparse.
In the case of underdetermined linear systems, our work can be related to the row action methods for compressed sensing presented in \cite{SraTropp:06}. This work differs from our approach in that it uses sequential projections onto the rows of $A$ followed by a thresholding step. \vspace{-0.1in}

\section{Sparse Randomized Kaczmarz}\label{sec:SparseRandKacz}\vspace{-0.1in}
We are interested in the case where the solution to the system of equations $Ax = b$ is sparse with at most $k < n$ non-zero entries. In this case, we assume that the signal $x$ is supported on a set $S_0 \subset \{1,\dots n\}$. Therefore, the measurement vector can be expressed as 
$$b = Ax = A_{S_0}x_{S_0},$$ 
where $A_{S_0}$ is the submatrix of $A$ whose columns are indexed by the entires of the set $S_0$, and $x_{S_0}$ is the restriction of $x$ to the set $S_0$.

Under the assumption that $S_0$ is known, it is more efficient to use the RK algorithm to recover $x_{S_0}$ from the system of equations $b = A_{S_0}x_{S_0}$, since the expected convergence rate then depends on the scaled condition number of $A_{S_0}$. Unfortunately, such an assumption is rarely (if ever) valid. However, it may be possible to estimate the set $S_0$ from the intermediate iterates of the RK algorithm. We propose to use the index set $S$ of the largest in magnitude $n-j+1$ entries of $x_{j-1}$ as the support estimate in iteration $j$. Our objective is then to orthogonally project the iterate $x_{j-1}$ onto the reduced space of $A_{S}x_{S} = b$.

To formulate this reduced projection, recall that $b(i) = \langle a_i, x\rangle$, where $x$ is the signal to be found. The Kaczmarz iteration can then be seen as applying the following split projection
\begin{equation}\label{eq:KaczProj}
\begin{array}{ll}
	x_j &= P_{a_i}(x) + P_{a_i^{\perp}}(x_{j-1})\\
		&= \frac{\langle a_i, x\rangle}{\langle a_i, a_i\rangle}a_{i}^T + \left(x_{j-1} - \frac{\langle a_i, x_{j-1}\rangle}{\langle a_i, a_i\rangle}a_{i}^T \right),
\end{array}
\end{equation}  
where $P_{a_i}(\cdot) = \frac{\langle a_i, \cdot\rangle}{\langle a_i, a_i\rangle}a_i^T$ is the projection operator onto $a_i$, and $P_{a_i^{\perp}}(\cdot)$ is the projection onto the orthogonal complement of $a_i$. 

Let $a_{i S}$ be the restriction of the row vector $a_i$ to the set $S$. Our objective is to apply the following split projection in each iteration
\begin{equation}\label{eq:IdealModKaczProj}
	x_j = P_{a_{i S}}(x) + P_{a_{i S}^{\perp}}(x_{j-1}),
\end{equation}
where $S$ is an estimate of the true support $S_0$. Such a projection step is only feasible if $S_0 \subseteq S$. However, it is possible that $S$ misses entries that are in $S_0$ in which case $b(i) = \langle a_{i S_0}, x\rangle \neq \langle a_{i S}, x\rangle$. Therefore, our algorithm uses the \emph{approximate} projection step
\begin{equation}
	x_j = \frac{\langle a_i, x\rangle}{\langle a_{i S}, a_{i S}\rangle}a_{i S}^T + P_{a_{i S}^{\perp}}(x_{j-1}),
\end{equation}
where $S = S_j$ is initialized to be the full row vector $\{1,\dots n\}$ and then gradually reduced by removing the indices that correspond to the entries of $x_{j-1}$ that are smallest in magnitude. Moreover, in order to allow entries in $S_0$ that may still be missing from $S$, we replace $a_{i S}$ with a weighted row vector $(\mathrm{w}_j\odot a_i)$, where the weight vector $\mathrm{w}_j$ is given by
$$\mathrm{w}_j(l) = \left\{
\begin{array}{ll} 1&, l \in S \\ \frac{1}{\sqrt{j}}&, l \in S^c\end{array}\right.$$
where $j$ is the iteration number. The operator $\odot$ denotes the Hadamard (element-wise) product.

The use of nonzero weight values on the set $S^c$ allows erroneous entries in $S$ to be removed and missing entries to be included in subsequent iterations. The modified sparse Kaczmarz iteration is then given by
\begin{equation}\label{eq:ModKaczProj}
\begin{array}{ll}
	x_j &= \frac{\langle a_i, x\rangle}{\langle \mathrm{w}_j\odot a_i, \mathrm{w}_j\odot a_i\rangle}(\mathrm{w}_j\odot a_i)^T + P_{(\mathrm{w}_j\odot a_i)^{\perp}}(x_{j-1}) \\ \\
	      &= x_{j-1} + \frac{b(i) - \langle \mathrm{w}_j\odot a_i, x_{j-1}\rangle}{\|\mathrm{w}_j\odot a_i\|_2^2} (\mathrm{w}_j\odot a_i)^T.
\end{array}
\end{equation}
Note that as $j \rightarrow \infty$, the weighted row vector $\mathrm{w}_j\odot a_i \rightarrow a_{i S}$. The full details of the sparse randomized Kaczmarz algorithm are summarized in Algorithm \ref{alg:SparseRK}.

\begin{algorithm}\caption{Sparse Randomized Kaczmarz algorithm}
\begin{algorithmic}[1]\label{alg:SparseRK}
\STATE \textbf{Input} $b = Ax$, the $m \times n$ matrix $A$, support estimate size $\hat{k}$, maximum iterations $J$
\STATE \textbf{Output} $x_j$
\STATE \textbf{Initialize} $S = \{1,\dots n\}$, $j = 0$, $x_0 = 0$  
\WHILE{$j \leq J$}
\STATE $j = j + 1$
\STATE Choose the row vector $a_i$ indexed by $i \in \{1,\dots m\}$ with probability $\frac{\|a_i\|_2^2}{\|A\|_F^2}$
\STATE Identify the support estimate $S$, such that $$S = \text{supp}(x_{j-1}|_{\max\{\hat{k}, n - j+1\}})$$
\STATE Generate the weight vector $\mathrm{w}_j$ such that $$\mathrm{w}_j(l) = \left\{
\begin{array}{ll} 1&, l \in S \\ \frac{1}{\sqrt{j}}&, l \in S^c\end{array}\right.$$
\STATE {\large$x_j = x_{j-1} + \frac{b(i) - \langle \mathrm{w}_j\odot a_i, x_{j-1}\rangle}{\|\mathrm{w}_j\odot a_i\|_2^2} (\mathrm{w}_j\odot a_i)^T$}
\ENDWHILE
\end{algorithmic}
\end{algorithm}
\vspace{-0.2in}

\section{Numerical Experiments}\label{sec:NumExp}\vspace{-0.1in}
We tested our sparse randomized Kaczmarz (SRK) algorithm by comparing its performance in recovering sparse solutions to the system $Ax = b$ with that of randomized Kaczmarz (RK) in the overdetermined system case, and with RK and standard $\ell_1$ minimization in the underdetermined system case. \vspace{-0.3in}

\subsection{Overdetermined systems}
\begin{figure*}[ht]
	\centering
	\includegraphics[width = 0.35\textwidth]{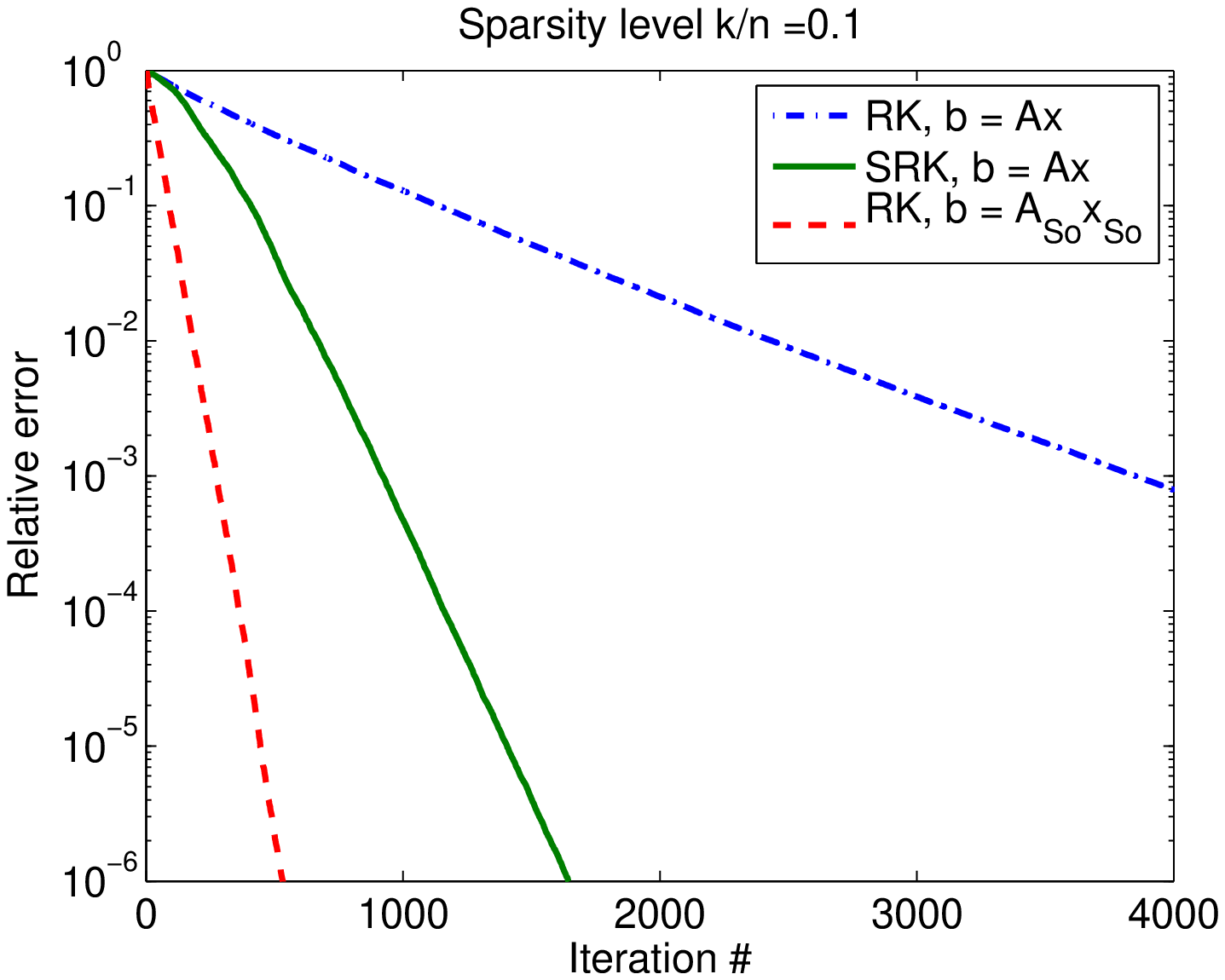}
	\includegraphics[width = 0.35\textwidth]{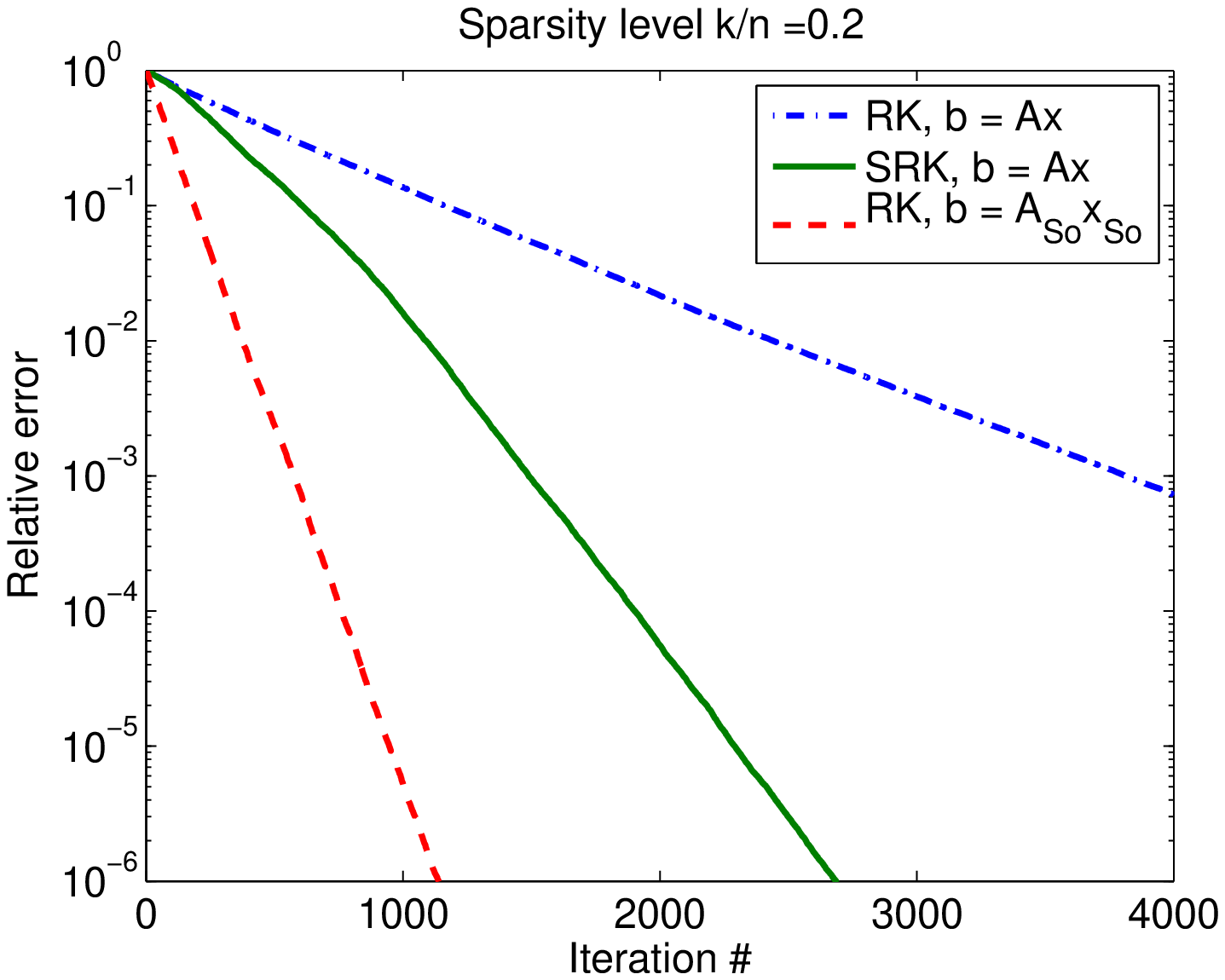}
	\includegraphics[width = 0.35\textwidth]{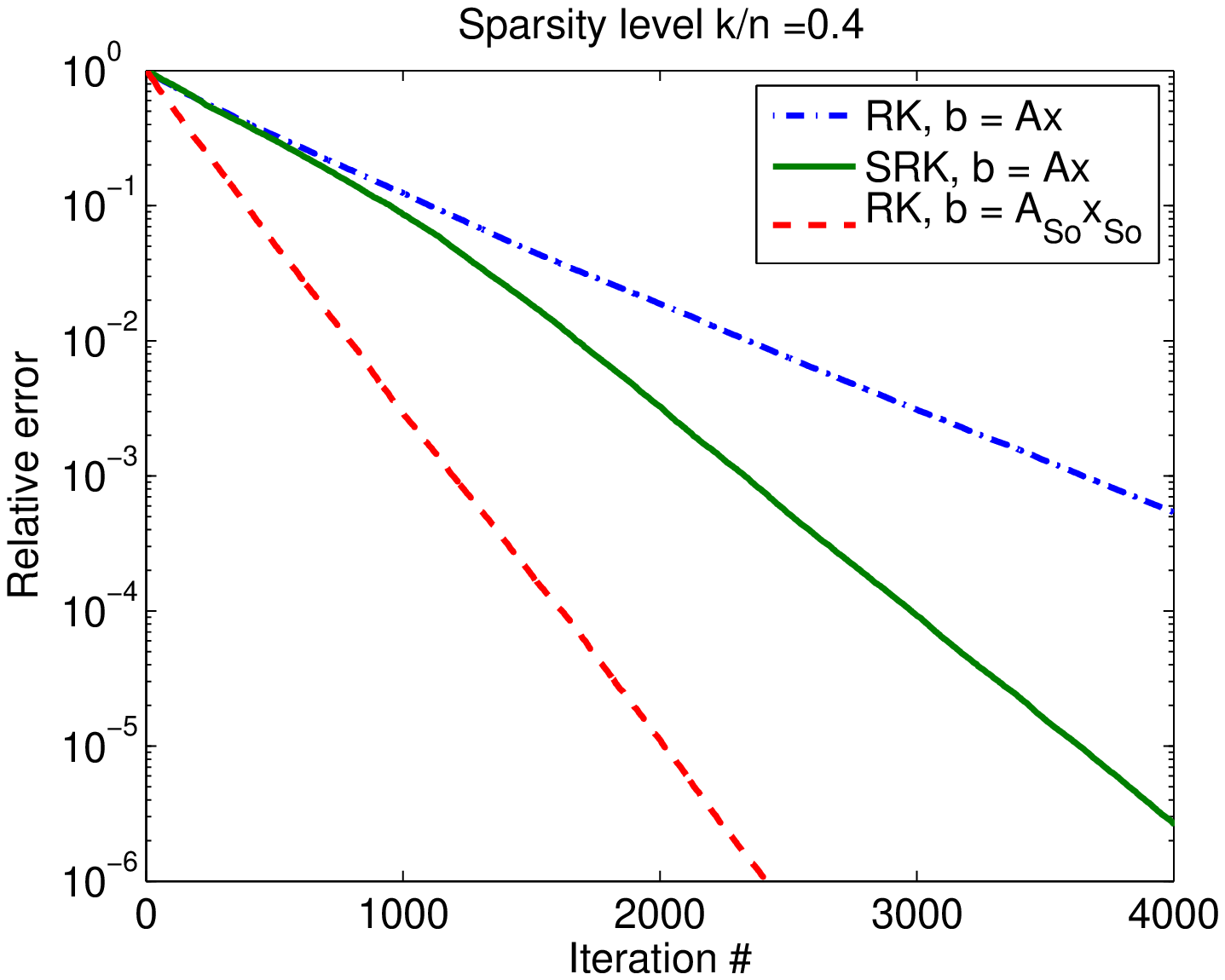}
	\includegraphics[width = 0.35\textwidth]{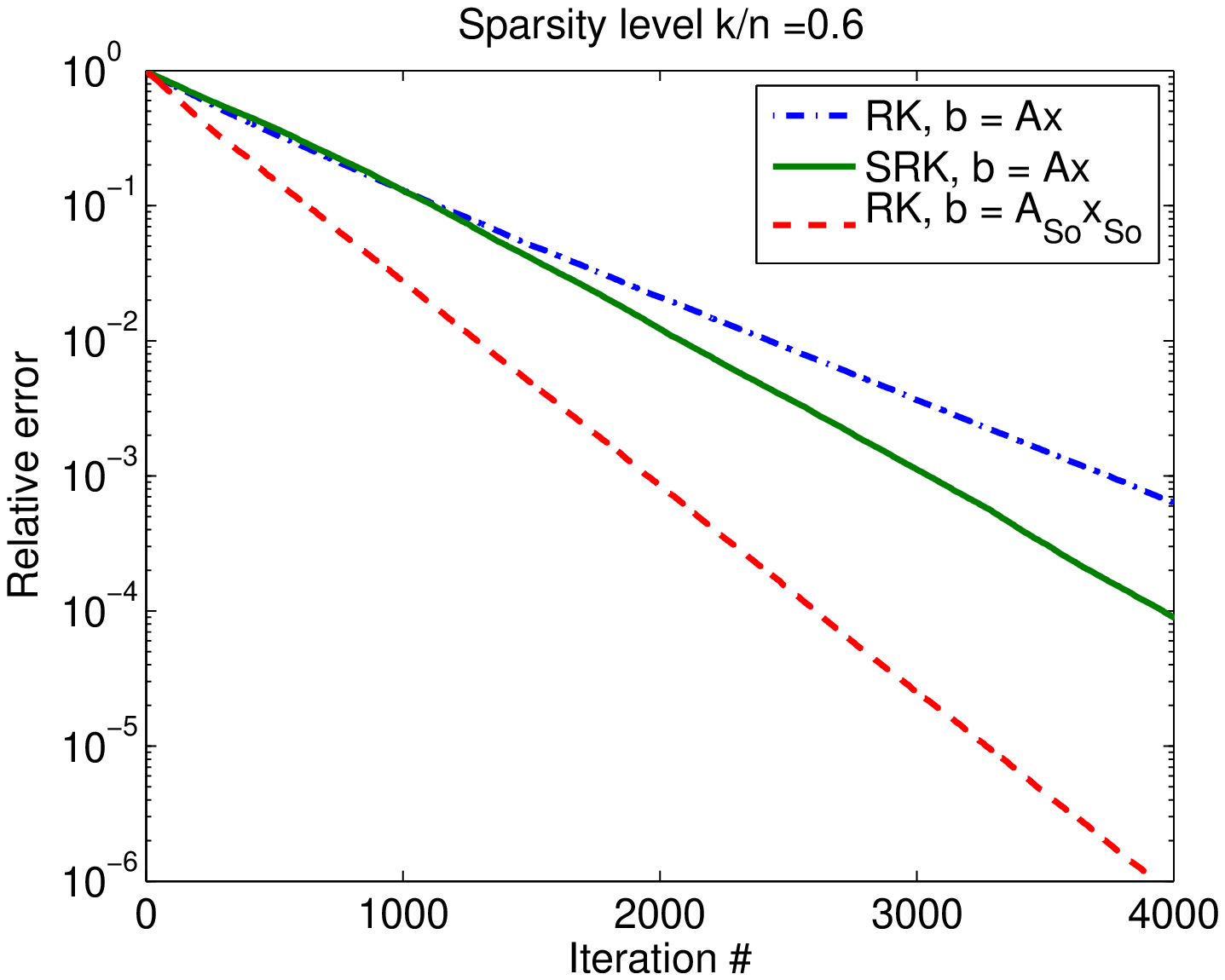}
	\caption{Average convergence rates over 20 runs of the proposed SRK algorithm and that of the RK algorithm for solving the overdetermined system $Ax = b$ with a sparse solution vector with sparsity $k/n$ for a $1000\times 200$ Guassian matrix $A$. The dashed red line shows the convergence rate of RK for the system $A_{S_0}x_{S_0} = b$, where $S_0$ is the support of the vector $x$.}\label{fig:Sparse_over}
\end{figure*}

\vspace{-0.1in}
We generate an $m \times n$ matrix $A$ with independent identically distributed (i.i.d.) Gaussian random entries, with $m = 1000$ and $n = 200$. We generate sparse signals $x$ with sparsity level $\frac{k}{n} \in \{0.1, 0.2, 0.4, 0.6\}$ and run each algorithm 20 times to recover the signal $x$. Figure~\ref{fig:Sparse_over} compares the convergence rates of SRK and RK averaged over the 20 experiments. We also show the convergence rate of RK applied to the reduced system $A_{S_0}x_{S_0} = b$ to illustrate the upper limit on the performance of SRK, where $S_0$ is the true support of $x$. The plots show that the SRK algorithm enjoys a significantly faster convergence rate than the RK algorithm, especially when the solution is very sparse.
\vspace{-0.1in}

\subsection{Underdetermined systems}
\begin{figure*}[hb]
	\centering
	\includegraphics[width = 0.33\textwidth]{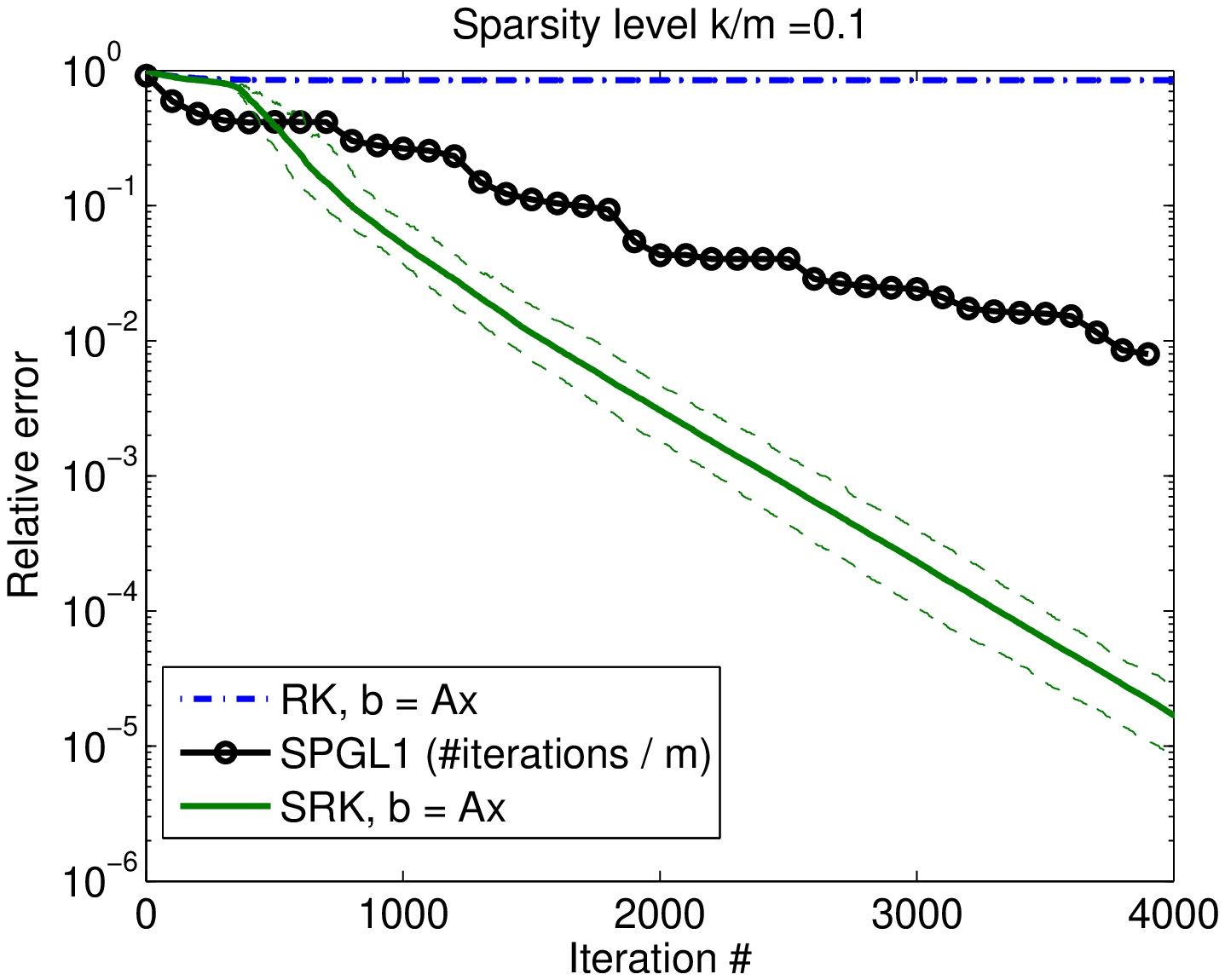}
	\includegraphics[width = 0.33\textwidth]{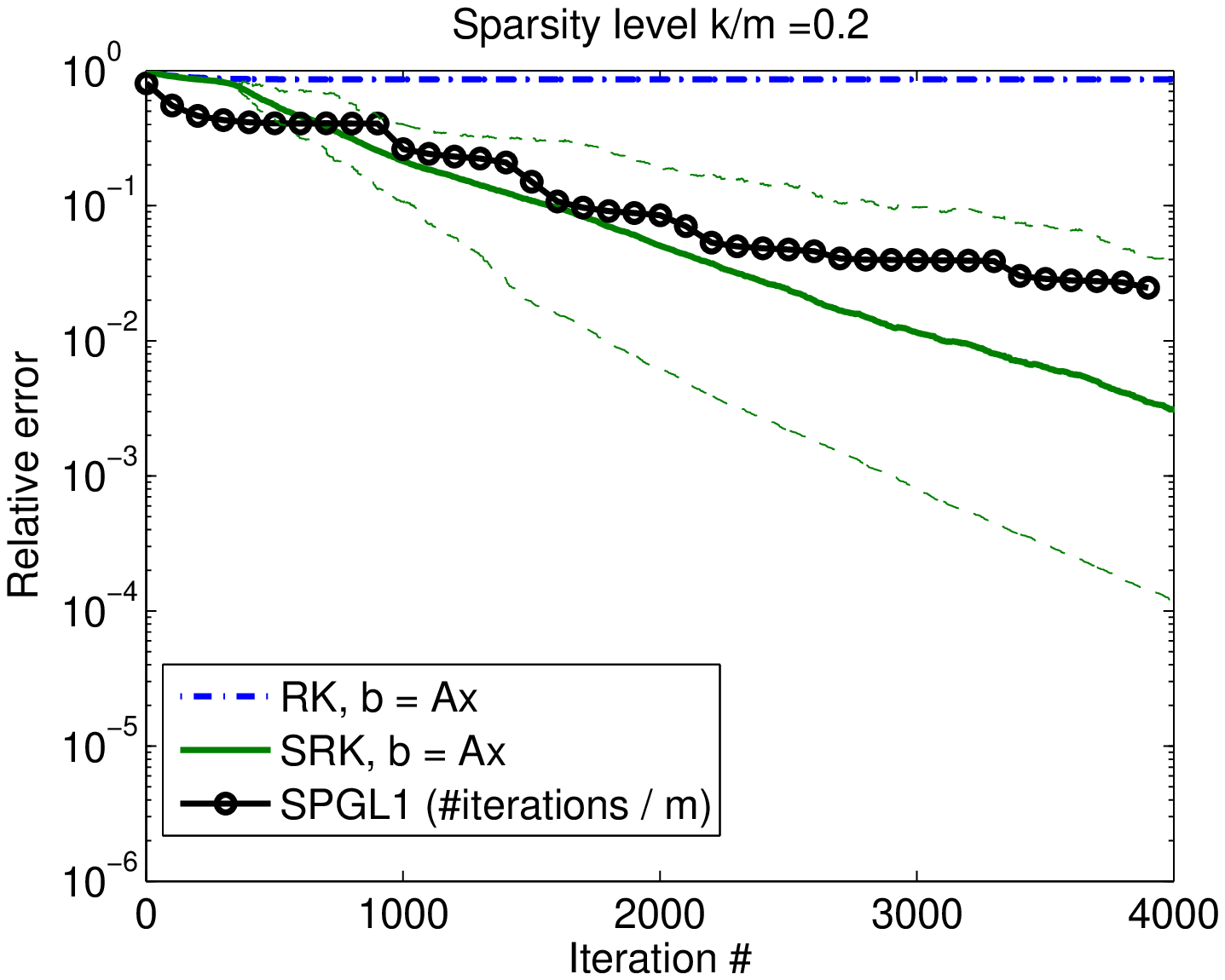}
	\includegraphics[width = 0.33\textwidth]{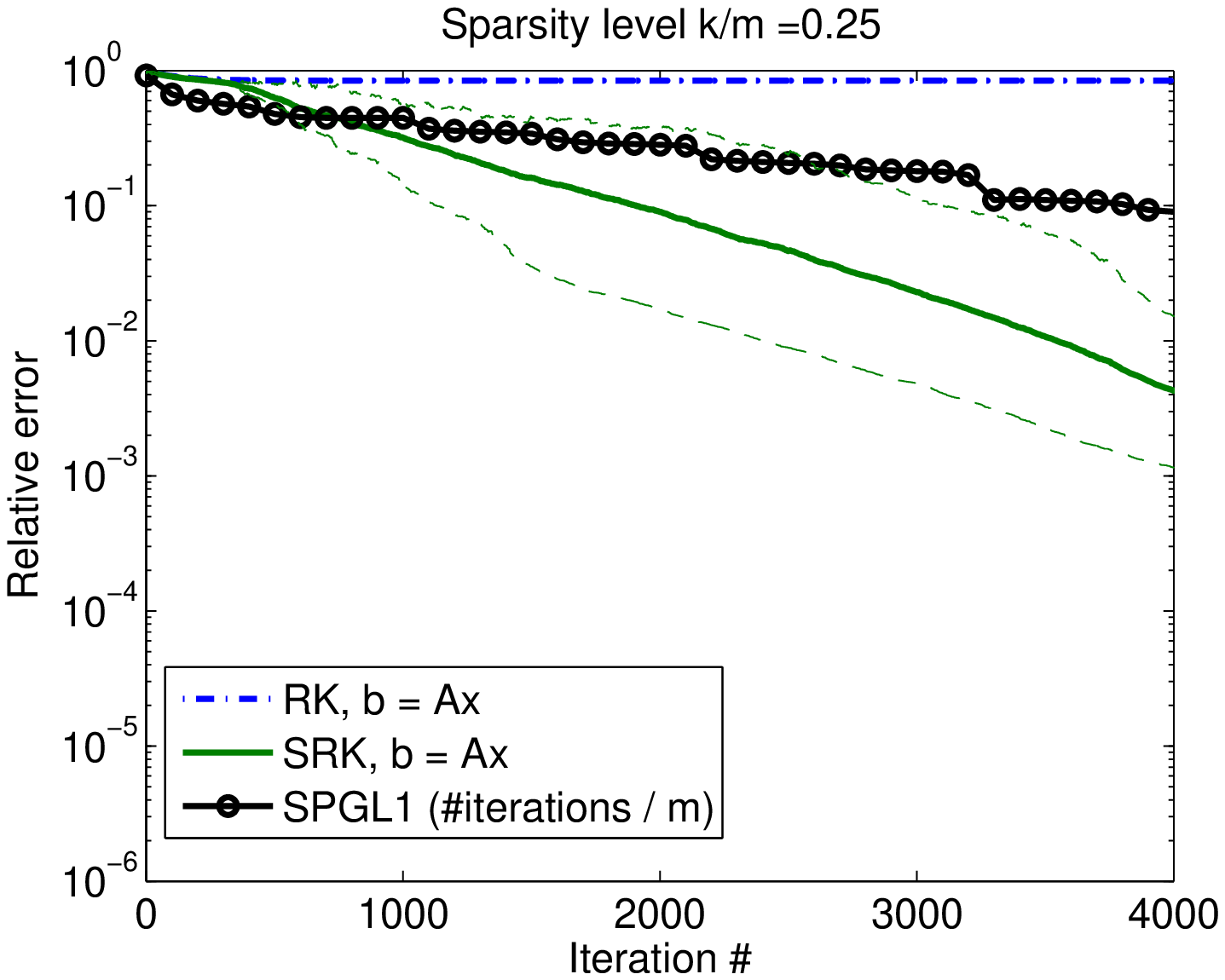}
	\caption{Average convergence rates over 20 runs of the proposed SRK algorithm and that of the RK algorithm for solving the underdetermined system $Ax = b$ with a sparse solution vector with sparsity $k/m$ for a $100\times 400$ Gaussian matrix $A$. The convergence rate of SPGL1 is shown where the number of displayed iterations is divided by $m$. The green dotted lines show the slowest and fastest convergence rates among the 20 runs.}\label{fig:Sparse_under}
\end{figure*}

\vspace{-0.1in}
In the underdetermined case, we generate a Gaussian matrix $A$ with dimensions $m = 100$ and $n = 400$. We also generate sparse signals $x$ with the sparsity level now determined with respect to the number of measurements $\frac{k}{m} \in \{0.1, 0.2, 0.25\}$ and run each algorithm 20 times to recover the signal $x$. The support estimate size of SRK is chosen as $\hat{k} = 0.25k, 0.2k,$ and $0.15k$, respectively. We compare the recovery result of SRK with the RK algorithm and with the spectral projected gradient for $\ell_1$ minimization (SPGL1) algorithm \cite{BergFriedlander:2008}. When comparing with the SPGL1 algorithm, we limit the number of SPGL1 iterations to the number of iterations of SRK divided by the number of measurements $m$. This allows us to keep the compare both algorithms at a reasonably equal number of vector-vector products. Figure~\ref{fig:Sparse_under} shows the convergence rates of SRK, RK, and SPGL1 under the fairness condition described above. It can be seen that the SRK algorithm can successfully recover sparse vectors from underdetermined systems within the same sparse recovery capabilities of $\ell_1$ minimization. However, we observe also that as the sparsity level increases, e.g. $k/m = 0.25$, the probability at which SRK succeeds in solving the problem decreases while it remains the same for SPGL1. The sparse recovery performance of SRK demonstrates its suitability for compressed sensing type problems. \vspace{-0.3in}

%
%

\small
\bibliographystyle{IEEEbib}
\bibliography{sparse}

\end{document}